\documentclass[12pt,a4paper]{article}
\usepackage[english]{babel}
\usepackage{amsmath}
\usepackage{amstext}
\usepackage{amssymb}
\usepackage{amsthm}
\usepackage{amsfonts}
\usepackage{fullpage}
\usepackage{graphicx}
\usepackage{epsf}
\usepackage[dvips]{epsfig}
\usepackage{stmaryrd}
\usepackage[T1]{fontenc}
\usepackage{latexsym}
\usepackage{psfrag}
\usepackage{cyrillic}

\title{On the finite blocking property}
\author{Thierry Monteil \footnote{Institut de Math\'ematiques de Luminy, CNRS UPR 9016,
Case 907, 163 Avenue de Luminy, 13288 Marseille cedex 09, France
-- E-Mail: monteil@iml.univ-mrs.fr -- Tel: +33 4 91 26 96 77 --
Fax : +33 4 91 26 96 55 --}}
\date{}

\newcommand{\N}{\mathbb N}
\newcommand{\R}{\mathbb R}
\newcommand{\C}{\mathbb C}
\newcommand{\Z}{\mathbb Z}
\newcommand{\Q}{\mathbb Q}
\newcommand{\defi}{\stackrel{\mbox {\tiny def}}{=}}

\newenvironment{appli}{\left( \begin{array}{ccc}}{\end{array} \right)}
\newcommand{\dans}{& \longrightarrow &}
\newcommand{\donne}{& \longmapsto &}
\newcommand{\ba}{\begin{appli}}
\newcommand{\ea}{\end{appli}}
\newtheorem{lem}{Lemma}[]
\newtheorem{prop}{Proposition}[]

\newtheorem{thm}{Theorem}[]

\newcommand{\B}{\mathbb B} 
\renewcommand{\H}{\mathcal H}

\newcommand{\T}{{\mathcal{T}}}
\newcommand{\extr}{\ \uparrow \!\! (\N,\N)}

\newcommand{\Int}{\stackrel{\circ}}

\newcommand{\card}{\mbox{card}}

\newcommand{\cge}{\xrightarrow[ n \rightarrow \infty]{}}

\renewcommand{\P}{{\mathcal{P}}}
\renewcommand{\S}{{\mathcal{S}}}
\newcommand{\E}{{\mathcal{E}}}
\renewcommand{\C}{{\mathcal{C}}}
\renewcommand{\B}{{\mathcal{B}}}

\renewcommand{\T}{{\mathcal{T}}}

\begin{document}

\maketitle


\begin{abstract}
\noindent A planar polygonal billiard $\P$ is said to have the
finite blocking property if for every pair  $(O,A)$ of points in
$\P$ there exists a finite number of ``blocking'' points $B_1,
\dots , B_n$ such that every billiard trajectory from $O$ to $A$
meets one of the $B_i$'s. Generalizing our construction of a
counter-example to a theorem of Hiemer and Snurnikov (see
\cite{Mo}), we show that the only regular polygons that have the
finite blocking property are the square, the equilateral triangle
and the hexagon. Then we extend this result to translation
surfaces. We prove that the only Veech surfaces with the finite
blocking property are the torus branched coverings. We also
provide a local sufficient condition for a translation surface to
fail the finite blocking property. This enables us to give a
complete classification for the L-shaped surfaces as well as to obtain 
a density result in the space of
translation surfaces in every genus~$g\geq 2$.\\


\end{abstract}

\clearpage

\section*{Introduction}

\begin{figure}[h!]
\begin{center}
\psfrag{p1=p2}{$p_1=p_2$} \psfrag{q1}{$\theta_1$}
\psfrag{q2}{$\theta_2$} \psfrag{q=0}{$\theta = 0$}
\includegraphics[width=.9 \linewidth]{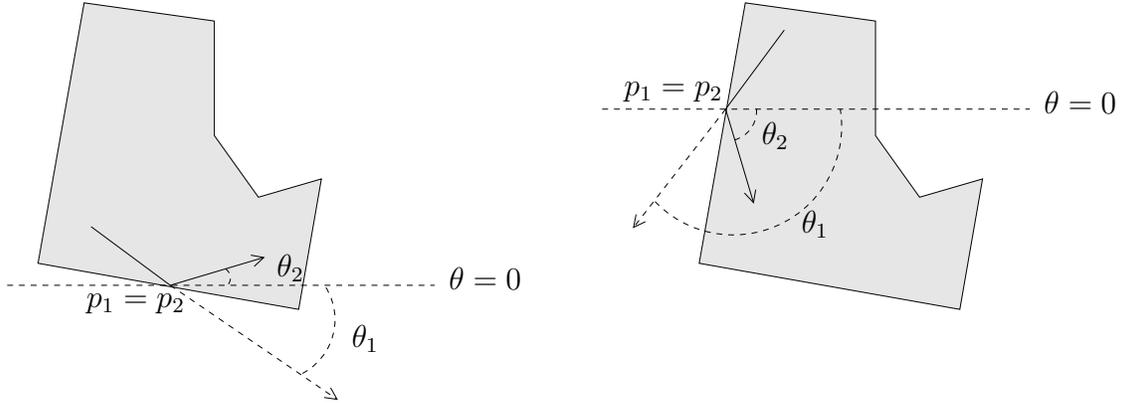}
\caption{\label{descartes} The Descartes law: \em{the incidence
angle equals the angle of reflection}.}
\end{center}
\end{figure}

When studying the motion of a point-mass in a polygonal billiard
$\P$, we work on the phase space $X=\P \times \mathbb{S}^1$
suitably quotiented: we identify the points $(p_1,\theta_1)$ and
$(p_2,\theta_2)$ if $p_1=p_2$ is on the boundary of $\P$ and if
the angles $\theta_1$ and $\theta_2$ are such that the Descartes
law of reflection is respected
(see Figure \ref{descartes}
).\\

The phase space enjoys also such a global decomposition in the
study of the dynamics on a translation surface.\\

There are essentially two points of view depending on whether the
variable is the first or the second projection:
\begin{enumerate}

\item We can fix one (or a finite number of) particular direction:
this corresponds for rational billiards to the study of the
directional flow in a translation surface (we are interested by
the ergodic properties depending on whether $\theta$ is a saddle
connection direction or not) (see \cite{KMS}, \cite{Ve},
\cite{MT}, \cite{Vo}). It is also useful for finding periodic
trajectories in irrational billiards by starting perpendicularly
to an edge (see \cite{ST}). This point of view is the most
studied.

\item We can also fix one (or a finite number of) point in $\P$
and look at which points we can reach when we let $\theta$ move.
This class of problems is called ``illumination problems''. The
first published question seems to appear in \cite{Kl} (see
\cite{KW} for a more precise story). The first published result in
this direction seems to be the paper of George W. Tokarski (see
\cite{To}) who finds a polygon that is not illuminable from every
point. Independently, Michael Boshernitzan constructed such an
example in a correspondence with Howard Masur (see \cite{Bos}).

\end{enumerate}

We are interested here in an illumination problem called the finite blocking property.\\

A planar polygon (resp. translation surface) $\P$ is said to have
the {\em finite blocking property} if for every pair  $(O,A)$ of points in $\P$,
there exists a finite number of 
points $B_1, \dots , B_n$  (different from $O$ and $A$) such that
every billiard trajectory (resp. geodesic) from $O$ to $A$ meets
one of the $B_i$'s.\\ \\

In this paper we will primarily focus on
translation surfaces. The paper is organized as follows: in
section \ref{background}, we will give some definitions and prove
that the finite blocking property is stable under branched
covering, stable under the Zemljakov-Katok's construction, and
stable under the action of $GL(2,\R)$. Section \ref{rqHS} is
devoted to the study of Hiemer and Snurnikov's proof, leading to
some comments about the finite property in the torus $\R^2 /
\Z^2$. In section \ref{local}, we prove a local sufficient
condition for a translation surface to fail the finite blocking
property (Lemma \ref{lem2cyln}). The aim of the next two sections
is to prove the following theorems:\\

{\bf Theorem \ref{th-regular}.} {\em Let $n\geq 3$ be an integer.
The following assertions are equivalent:
\begin{itemize}
\item the regular $n$-gon has the finite blocking property. \item
the right-angled triangle with an angle equal to  $\pi/n$ has the
finite blocking property. \item $n\in\{3,4,6\}$.
\end{itemize}}

{\bf Theorem \ref{th-Veech}.} {\em A Veech surface has the finite
blocking property if and only if it is a torus covering, branched
at only one point.}\\



The last section is devoted to some other applications of Lemma
\ref{lem2cyln}, such as\\

{\bf Proposition \ref{L-prop}.} {\em Let $a$ and $b$ be two
positive real numbers. Then the L-shaped surface $L(a,b)$ has the
finite blocking property if and only if $(a,b)\in \Q^2$.}\\

{\bf Theorem \ref{th-density}.} {\em In genus $g\geq 2$, the set
of translation surfaces that fail the finite blocking property is
dense in every stratum.}

\vspace{2cm}


{\bf Acknowledgement:} I would like to thank Martin Schmoll for
introducing me to the subject, Kostya Kokhas and Serge Troubetzkoy
for historical comments, Anton Zorich  for helpful discussions and 
Pascal Hubert for encouragements to write this paper.

\clearpage

\section{Definitions and first results} \label{background}
\subsection{Translation surfaces and geodesics}

A {\em translation surface} is a triple $(\S, \Sigma, \omega)$
such that $\S$ is a topological compact connected surface,
$\Sigma$ is a finite subset of $\S$ (whose elements are called
{\em singularities}) and $\omega = (U_i,\phi_i)_{i\in I}$ is an
atlas of $\S \setminus \Sigma$ (consistent with the topological
structure on $\S$) such that the transition maps (i.e. the $\phi_j
\circ \phi_i^{-1} : \phi_i(U_i\cap U_j) \rightarrow \phi_j(U_i\cap
U_j)$ for $(i,j)\in I^2$) are translations. This atlas gives to
$\S \setminus \Sigma$ a Riemannian structure; we therefore have
notions of length, angle, measure, geodesic... We assume moreover
that $\S$ is the completion of $\S \setminus \Sigma$ for this
metric.
We will sometimes use the notation $(\S, \Sigma)$ or simply $\S$
to refer to $(\S, \Sigma, \omega)$. A singularity $\sigma \in
\Sigma$ is said to be {\em removable} if there exists an atlas
$\omega' \supset \omega$ such that $(\S, \Sigma \setminus \{
\sigma \}, \omega')$ is a translation
surface.\\

There are many conventions about what happens if a geodesic
$\gamma$ meets one singularity $\sigma \in \Sigma$. Some people
want to stop $\gamma$ here; some other people want to extend
$\gamma$ with a multi-geodesic path; other people want to extend
$\gamma$ after $\sigma$ if and only if $\sigma$ is an
removable singularity. 
{\em The finite blocking property does not depend on the
convention} since, for every pair  $(O,A)$ of points in $\S$, we
can always add the set $\Sigma \setminus \{O,A\}$ (that is finite)
to a blocking configuration. Therefore, if $(S,\Sigma)$ is a
translation surface, if $\E$ is a finite subset of $\S$ that we
want to add in $\Sigma$ as removable singularities, then
$(S,\Sigma)$ has the finite blocking property if and only if
$(S,\Sigma \cup \E)$ has it.

\subsection{Branched coverings} \label{covering}

A {\em branched covering} between two translation surfaces is a
mapping $\pi : (\S,\Sigma) \rightarrow (\S',\Sigma')$ that is
a topological branched covering that locally preserves the
translation structure.

\begin{prop}\label{prop-covering}
Let $\pi : (\S,\Sigma) \rightarrow (\S',\Sigma')$ be a covering of
translation surfaces branched on a finite set $\mathcal{R}'
\subset \S'$. Then $\S$ has the finite blocking property if and
only if $\S'$ has.
\end{prop}

{\bf Proof:} $\Rightarrow$ : Suppose that $\S$ has the finite
blocking property. Let $(O',A')$ be a pair of points in $\S'$. Let
$O$ be a point chosen in  $\pi^{-1}(\{O'\})$. If $A\in
\pi^{-1}(\{A'\})$, there exists a finite set $\B_A$ of points in
$\S \setminus \{O,A\}$ such that every geodesic in $\S$ from $O$
to $A$ meets $\B_A$. Let $$\B' \defi \left(  \bigcup_{A\in
\pi^{-1}(\{A'\})} \pi( \B_A)  \cup \mathcal{R}' \right) \setminus
\{O',A'\}.$$ Let $\gamma' : [a,b] \rightarrow \S'$ be a geodesic
from $O'$ to $A'$. Up to a restriction, we can suppose that
$\gamma'(]a,b[) \cap \{O',A'\} = \emptyset$. Suppose by
contradiction that $\gamma'([a,b]) \cap \B' = \emptyset$. In
particular, $\gamma'(]a,b[) \cap \mathcal{R}' = \emptyset$. So,
$\gamma'$ can be lifted to a geodesic $\gamma: [a,b] \ \rightarrow
\S'$ from $O$ to some $A\in \pi^{-1}(\{A'\})$ such that $\pi \circ
\gamma = \gamma'$. Then, there exists $t \in ]a,b[$ such that
$\gamma(t) \in \B_A$. Hence $\gamma'(t) \in \B'$, leading to a
contradiction. So, $\B'$ is a finite blocking configuration and
$\S'$ has the finite blocking property. \\

$\Leftarrow$ : Suppose that $\S'$ has the finite blocking
property. Let $(O,A)$ be a pair of points in $\S$. Let $O'\defi
\pi(O)$ and $A'\defi \pi(A)$. There exists a finite set $\B'
\subset \S' \setminus \{O',A'\}$ such that every geodesic in $\S'$
from $O'$ to $A'$ meets $\B'$. Let $$\B \defi \pi^{-1}(\B')
\subset \S \setminus \{O,A\}.$$
Let $\gamma : [a,b] \rightarrow \S$ be a geodesic from $O$ to $A$.
$\gamma$ can be pushed to a geodesic $\gamma' \defi \pi \circ
\gamma: [a,b] \ \rightarrow \S'$ from $O'$ to $A'$. Then, there
exists $t \in ]a,b[$ such that $\gamma'(t) \in \B'$. Hence
$\gamma(t) \in \B$. So, $\B$ is a finite blocking configuration
and $\S$ has the finite blocking property.

\hfill $\square$

\subsection{Rational billiards {\em vs} translation surfaces}

Let $\P$ denote a polygon in $\mathbb{R}^2$, whose set of vertices
is denoted by $V$. Let $\Gamma \subset O(2,\R)$ be the group
generated by the linear parts of the reflections in the sides of
$\P$. When $\Gamma$ is finite, we say that $\P$ is a {\em rational
polygonal billiard}. When $\P$ is simply connected, $\P$ is
rational if and only if all the angles between edges are rational
multiples of $\pi$.\\

A classical construction due to Zemljakov and Katok (see
\cite{ZK}, \cite{MT}) allows us to associate to each rational
billiard $\P$ a translation surface $ZK(\P)$ as follows:

Let $(P_{\gamma})_{\gamma \in \Gamma}$ be a family of $|\Gamma|$
disjoint copies of $\P$, each $P_{\gamma}=\gamma(\P)$ being
rotated by the element $\gamma \in \Gamma$. If $\gamma\in\Gamma$,
if $e$ is an edge of $P_{\gamma}$, let $\delta\in\Gamma$ be the
linear part of the reflection in $e$; we identify $e\in
P_{\gamma}$ with $\delta(e)\in P_{\delta \gamma}$. We set:
$$
ZK(\P) \defi   \bigsqcup_{\gamma \in \Gamma} P_{\gamma} \
\mbox{\bf \Large{/}}  \sim
$$
where $\sim$ is the relation above. The translation structure of
each $\Int{P_{\gamma}} \in \R^2$ can be extended to an atlas of
$\bigcup_{\gamma \in \Gamma} P_{\gamma} \setminus \gamma(V)$, that
gives to $ZK(\P)$ a translation structure whose set of
singularities is $\Sigma=\bigcup_{\gamma \in \Gamma} \gamma(V)$.\\

In other terms, $\P$ tiles $ZK(\P)$ which can be written as
$ZK(\P) = \bigcup_{\gamma\in\Gamma} \psi_{\gamma}(\P)$ where the
$\psi_{\gamma}$'s are isometries. Let $$\displaystyle \psi \defi
\ba  ZK(\P) \dans \P \\ x \donne
({\psi_{\gamma}}_{|_{\psi_{\gamma}(\P)}} )^{-1} (x) \ \ \mbox{ if
} x\in \psi_{\gamma}(\P) \ea .$$ $\psi$ is well defined since if
$x\in \psi_{\gamma}(\P) \cap \psi_{\delta}(\P)$, then
$({\psi_{\gamma}}_{|_{\psi_{\gamma}(\P)}} )^{-1} (x) =
({\psi_{\delta}}_{|_{\psi_{\delta}(\P)}} )^{-1} (x)$ (this is just
the compatibility with $\sim$).
Moreover, $\psi$ is a piecewise isometry.\\

\begin{prop}\label{prop-ZK}
Let $\P$ be a rational polygonal billiard. Then $ZK(\P)$ has the
finite blocking property if and only if $\P$ has.
\end{prop}

{\bf Proof:} it is very similar to the proof given in subsection
\ref{covering} ($\psi$ plays the role of $\pi$).

Indeed, for the direction $\Rightarrow$, if $(O',A')$ is a pair of
points in $\P$, if $O$ is chosen in $\psi^{-1}(\{O'\})$, then for
each $A$ in $\psi^{-1}(\{A'\})$, there exists a finite set $\B_A$
of points in $ZK(\P) \setminus \{O,A\}$ such that every geodesic
in $ZK(\P)$ from $O$ to $A$ meets $\B_A$.
Then $$\B' \defi \left(  \bigcup_{A\in \psi^{-1}(\{A'\})} \psi(
\B_A)  \cup V \right) \setminus \{O',A'\} =\left(
\bigcup_{\gamma\in\Gamma} \psi (\B_{\psi_{\gamma} (A')})  \cup V
\right) \setminus \{O',A'\}$$
is a finite blocking configuration between $O'$ and $A'$, thus $\P$ has the finite blocking property. \\

For the direction $\Leftarrow$, if $(O,A)$ is a pair of points in
$ZK(\P)$, there exists a finite set $\B'$ of points in $\P
\setminus \{\psi(O),\psi(A)\}$ such that every billiard path in
$\P$ from $\psi(O)$ to $\psi(A)$ meets $\B'$.
Then $$\B \defi \psi^{-1}(\B') \subset ZK(\P) \setminus \{O,A\}$$
is a finite blocking configuration between $O$ and $A$, thus
$ZK(\P)$ has the finite blocking property.

\hfill $\square$ \\

{\bf Unfolding a billiard table:} combining propositions
\ref{prop-covering} and \ref{prop-ZK}, we have (see \cite{MT}) the:

\begin{prop}\label{prop-depliage}
Let $\P$ and $\P'$ be two rational polygonal billiards such that
$\P$ is obtained by reflecting $\P'$ at its edges finitely many
times, without overlapping (we allow some barriers along parts of
some sides of copies of $\P'$ inside $\P$). Since there exists a
branched covering from $ZK(\P)$ to $ZK(\P')$, then $\P$ has the
finite blocking property if and only if $\P'$ has.\\
\end{prop}

\mathversion{bold}
\subsection{Action of $GL(2,\R)$} \label{GL2R}
\mathversion{normal}

If $A\in GL(2,\R)$, we can define the translation surface $$A .
(\S, \Sigma,(U_i,\phi_i)_{i\in I}) \defi (\S, \Sigma, (U_i,A \circ
\phi_i)_{i\in I});$$ hence we have an action of $GL(2,\R)$ on the
class of translation surfaces. We classically consider only
elements of $SL(2,\R)$ (see \cite{MT}), but we do not need to
preserve area here.

\begin{prop}\label{prop-GL2R}
Let $\S$ be a translation surface and $A$ be in $GL(2,\R)$. Then
$\S$ has the finite blocking property if and only if $A.\S$ has.
\end{prop}

{\bf Proof:} such an action sends geodesic to geodesic.

\hfill $\square$ \\

To summarize this section, we can say that the finite blocking
property enjoys many properties of {\em stability}. As an
illustration, it suffices to apply successively propositions
\ref{prop-ZK},~\ref{prop-covering},~\ref{prop-GL2R},~\ref{prop-ZK},~\ref{prop-depliage} 
and to follow a construction of
\cite{Mc}, to reduce the problem of the finite blocking property
for the billiard in the regular pentagon to the problem in the
$\bot$-shaped billiard studied in \cite{Mo} with parameters
$\alpha=1+2\cos (2\pi/5)$, $L_1=1$, $L_2=2\cos (2\pi/5)$. It is
proved there that this billiard table fails the finite blocking
property ($\alpha \in \R \setminus \Q$).

\clearpage

\section{Some remarks around Hiemer and Snurnikov's proof} \label{rqHS}
In their article \cite{HS}, Philipp Hiemer and Vadim Snurnikov
tried to
prove that any rational billiard $\P$ has the finite blocking property.\\

For this, they use the subgroup $G_\P$ of $Isom(\R^2)$ generated
by the reflections at the edges of the polygon $\P$. In the proof
of their theorem $5$, they construct a finite number of points in
$\R^2$ called the $P_{i,\lambda}$'s and {\em choose} a blocking
point arbitrarily in each orbit of the $P_{i,\lambda}$'s under the
action of the group $G_\P$ on $\R^2$ (the $(i,\lambda)$'s belong
to $G_\P / T_\P  \times  \{0,1/2\}^{\dim_{\Q} T_\P}$).\\

\begin{figure}[h]
\begin{center}
\psfrag{a}{$\sqrt{2}$} \psfrag{1}{$1$}
\includegraphics[scale=0.7]{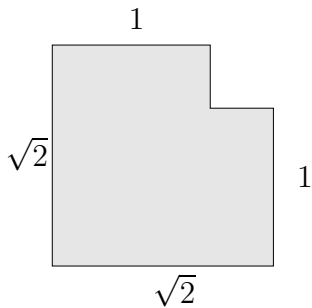}
\caption{\label{gpe-dense} A rational polygonal billiard whose
translation group is dense in $\R^2$.}
\end{center}
\end{figure}

The polygon $\P$ drawn on Figure \ref{gpe-dense} is a rational one
but the subgroup $T_\P$ of $G_\P$ consisting of translations is
dense in $\R^2$ (identified with the group of all translations of
$\R^2$), since $\Z+\sqrt{2}\Z$ is dense in $\R$. Hence the orbit
of any point of the plane under $G_\P$ is dense in the
plane and therefore in $\P$ (which is the closure of an open set).\\

\begin{figure}[h]
\begin{center}
\psfrag{a}{$\sqrt{2}$} \psfrag{1}{$1$} \psfrag{U}{$U$}
\psfrag{A}{$A$} \psfrag{O}{$O$}
\includegraphics[scale=0.7]{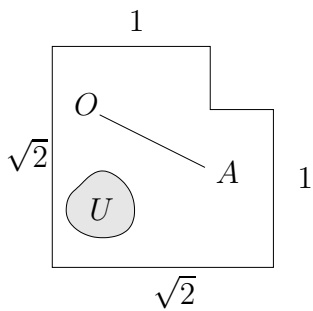}
\caption{\label{bug-demo} The direct path from $O$ to $A$ does not
meet the points in $U$.}
\end{center}
\end{figure}

Now, if we take two points $O$ and $A$ in $\P$ such that the
segment $[O,A]$ is included in $\P$, we can {\em choose} all the
blocking points in a small open set $U \subset \P$ which does not
intersect $[O,A]$. Such points cannot block the direct path from $O$ to
$A$ (see Figure \ref{bug-demo}).\\

Hence the proof in \cite{HS} does not work. In fact, we have shown
in \cite{Mo} that the billiard table drawn in Figure
\ref{gpe-dense} fails the finite blocking property. Meanwhile, the
proof given in \cite{HS} (theorem $5$) works for rational polygons
$\P$ such that $T_\P$ is discrete (hence lattice) (such billiards
are called {\em almost integrable}). Indeed, instead of choosing
one point in each orbit of the $P_{i,\lambda}$'s in $\P$, it
suffices to take {\em all} the points of the orbits of the
$P_{i,\lambda}$'s that lie in $\P$ (and that are distinct from $O$
and $A$).

We can bound (badly but uniformly) the number of such points as follows: \\

\begin{prop} \label{bound-almost}
If $\P$ is an almost integrable rational polygonal billiard with
angles of the form $\pi \frac{p_i}{q}$ ($1\leq i\leq n$), if the
diameter of $\P$ is $D$, if $(v_1,v_2)$ is a basis of $T_\P$ (that
is a lattice), then for every pair $(O,A)$ in $\P$, there exist a
set of at most $\displaystyle 8q \left[ \frac{\pi
(\frac{D}{\sqrt{3}} + \frac{\|v_1\| + \|v_2\|}{2})^2
}{\det(v_1,v_2)}  \right] $ blocking points.

\end{prop}

{\bf Proof:} include $\P$ in a disk of radius
$\frac{D}{\sqrt{3}}$, enlarge this disk by $\frac{\|v_1\| +
\|v_2\|}{2}$ and look at the area
($8q$ is larger than $\card (G_\P / T_\P \times
\{0,1/2\}^{\{v_1,v_2\}})$).

\hfill $\square$ \\

In fact, the conditions of this proposition (the fact of being an
almost integrable billiard) work only for the polygons constructed
by reflecting $\C$ at its edges finitely many times (we can add
barriers along interior edges), where $\C$ is one of the following
elementary polygons:

\begin{itemize}
  \item the right-angled isosceles triangle with angles $(\pi/2 , \pi/4 , \pi/4)$
  \item the half-equilateral triangle with angles $(\pi/2 , \pi/3 , \pi/6)$
  \item any rectangle
\end{itemize}

This fact is a direct consequence of the theory of reflection
groups and chamber systems (see \cite{Bou} for
an extensive study or \cite{Car} for a brief introduction). The three
elementary polygons correspond to the types $\tilde{C}_2$,
$\tilde{G}_2$ and
$\tilde{A}_1 \times \lambda \tilde{A}_1$ ($\lambda > 0$).\\

For example, the equilateral triangle, the square and the regular
hexagon have the finite blocking property. \\

Note that the translation surface associated to the square is a
(flat) torus. We can notice that the translation surface
associated to an almost integrable polygon is a torus branched
covering (this is easy for one of the elementary polygons, the
rest follows by reflecting). As we have seen in subsection
\ref{covering}, the finite blocking property is preserved by
branched covering and, since $4$ points suffice to block every
geodesic between $2$ fixed points in the torus, we have another
bound for the number of blocking points
for such billiards that depends on the degree of the covering.\\

With this remark, we can see that Proposition \ref{bound-almost}
can be seen as a consequence of the only fact that the square
billiard (or equivalently the translation surface $\R^2 / \Z^2$)
has the finite blocking property. This last result seems to appear
for the first time in the Leningrad's Olympiad in 1989 selection
round, $9$th form (see \cite{Fo}). The problem was the following:
\begin{quote}
``Professor Smith is standing in the squared hall with mirror
walls. Professor Jones wants to place in the hall several students
in such a way that professor Smith could not see from his place
his own mirror images. Is it possible? (Both professors and
students are points, the students can be placed in corners and
walls).''
\end{quote}
The author of this problem was Dmitrij Fomin. None of the school
students solved it. The booklet of the olympiad contains an
answer: $16$ students, and an example of this arrangement in
coordinates.\\

Another consequence of the fact that the torus $\R^2 / \Z^2$ has
the finite blocking property is the 

\begin{prop} There exists a translation surface with the property
that every geodesic going from a singularity to itself has to meet
first another singularity. In other words, this surface does not
have any saddle connection going from a singularity to itself.
\end{prop}

\begin{figure}[h!]
\begin{center}
\psfrag{a}{$A$} \psfrag{b}{$B$} \psfrag{c}{$C$} \psfrag{d}{$D$}
\psfrag{s}{$\S$} \psfrag{t}{$\R^2 / \Z^2$} \psfrag{pi}{$\pi$}
\psfrag{00}{$(0,0)$} \psfrag{01}{$(0,1)$} \psfrag{10}{$(1,0)$}
\includegraphics[width=.9 \linewidth]{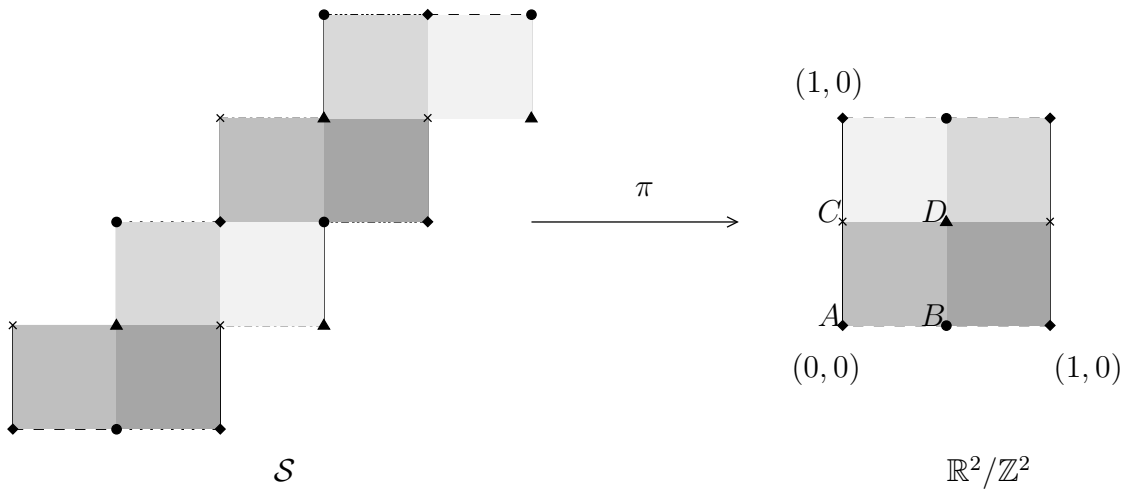}
\caption{\label{anton} Covering between $\S$ and $\R^2 / \Z^2$
(the opposite vertical lines are identified by horizontal
translation, the vertical lines with the same style are
identified).}
\end{center}
\end{figure}

{\bf Proof:} let $\S$ be the translation surface drawn in Figure
\ref{anton}.

The coloring  allows us to see how to construct a covering $\pi :
\S \rightarrow \R^2 / \Z^2$ of degree $2$, branched at the points
$A=(0,0)$, $B=(1/2,0)$, $C=(0,1/2)$, $D=(1/2,1/2)$. The
singularities of $\S$ are located at the preimages of those four points.\\

The finite blocking property in the square says that three points
in $\{A,B,C,D\}$ block every geodesic from the fourth point to
itself. Hence, by lifting, every geodesic from the singularity
$\pi^{-1}(A)$ to itself in $\S$ has to meet one of the other
singularities $\pi^{-1}(B)$, $\pi^{-1}(C)$ or $\pi^{-1}(D)$, and
by symmetry it works for the three other singularities.
\hfill $\square$ \\

\clearpage

\section{A local lemma}\label{local}

A {\em subcylinder} $\C$ is an isometric copy of $\R / w \Z \times
]0,h[$ in a translation surface $\S$ ($w>0$, $h>0$). $w$ and $h$
are unique and called the {\em width} and the {\em height} of
$\C$.

The images in $\S$ of $\R / w \Z \times \{0\}$ and  $\R / w \Z
\times \{h\}$ (which are well-defined if we extend the isometry,
which is uniformly continuous, by continuity in $\S$) are called
the {\em sides} of $\C$.

The {\em direction} of $\C$ is the direction of the image of $\R /
w \Z \times \{ h /2 \}$ (which is a closed geodesic).

A  {\em cylinder} is a maximal subcylinder (for the inclusion). By
maximality, each side of a cylinder must contain at least one
singularity and is a finite union of saddle connections.


Let $\C_1$ and $\C_2$ be two cylinders with width $w_1$ (resp.
$w_2$) and height $h_1$ (resp. $h_2$) in a translation surface
$\S$.

$\C_1$ and $\C_2$ are said to be {\em parallel} if their
directions are parallel.

$\C_1$ and $\C_2$ are said to be {\em commensurable} if the ratios
$w_1/h_1$ and $w_2/h_2$ are commensurable.\\

\begin{figure}[!h]
\begin{center}
\psfrag{w1}{$w_1$} \psfrag{h1}{$h_1$} \psfrag{w2}{$w_2$}
\psfrag{h2}{$h_2$} \psfrag{C1}{$\C_1$} \psfrag{C2}{$\C_2$}
\psfrag{O}{$(0,0)$} \psfrag{l}{$l$}
\includegraphics[scale = 0.6]{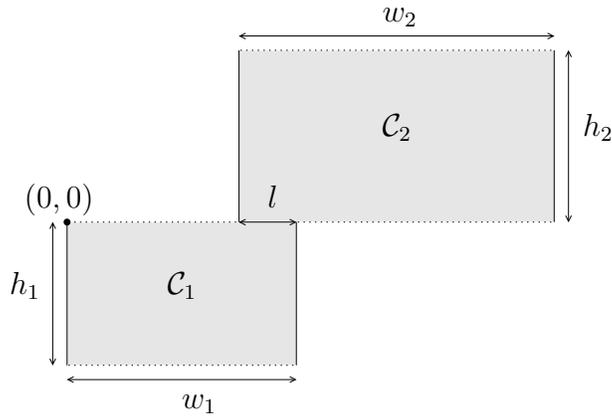}
\caption{\label{2cylinder.eps} The opposite vertical sides are
identified by translation while the dotted horizontal ones are
glued with the rest of the surface; we take $(0,0)$ as the origin;
$l>0$.}
\end{center}
\end{figure}

We will study the case where $\C_1$ and $\C_2$ are two different
parallel cylinders whose closures have a nontrivial (i.e. not reduced to a
finite set) intersection. The situation can be described as in
Figure \ref{2cylinder.eps}.

\begin{lem}\ \label{lem2cyln}
Let $\S$ be a translation surface that contains two different
parallel cylinders $\C_1$ and $\C_2$ whose closures have a nontrivial
intersection. If their widths are uncommensurable, then $\S$ fails
the finite blocking property.
\end{lem}

{\bf Proof:} Up to a vertical dilatation, we can assume that $h_1$
and $h_2$ are greater than $1$. In the system of coordinates given
in Figure \ref{2cylinder.eps}, we set $O=(w_1-l/2, -1)$ and
$A=(w_1-l/2,1)$.

Since $\frac{w_1}{w_2}$ is a positive irrational number, $\N^* -
\frac{w_1}{w_2} \N^*$ is dense in $\R$ so there exists two
positive integer sequences $(p_n)_{n \in \N}$ and $(q_n)_{n \in
\N}$ such that:

\begin{itemize}
  \item $q_n$ is strictly increasing
  \item $p_n w_2 - q_n w_1 \in ]-l, l[$
\end{itemize}

For $n \in \N$, let $\gamma_n$ be the geodesic starting form $O$
with slope
$$\frac{2}{q_n w_1  +  p_n w_2 }
=\frac{1}{ q_n w_1 + \lambda_n}=\frac{1}{ p_n w_2 - \lambda_n}$$
where $\lambda_n = (p_n w_2 - q_n w_1)/2 \in ]-l/2,l/2[$.

\begin{figure}[!h]
\begin{center}
\psfrag{totalslope}{\small $ q_n w_1  +  p_n w_2 $}
\psfrag{slope1}{\small $ q_n w_1 + \lambda_n $}
\psfrag{slope2}{\small $ p_n w_2 - \lambda_n $}
\includegraphics[width=.9 \linewidth]{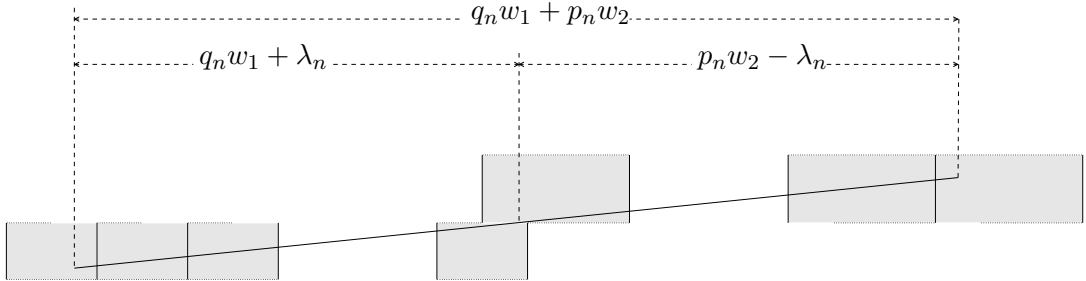}
\caption{\label{unfold-lemme.eps} The unfolding procedure.}
\end{center}
\end{figure}

So, we can check by unfolding the trajectory in the universal
cover of $\S$ (see Figure \ref{unfold-lemme.eps}) that $\gamma_n$
passes $q_n$ times through the line $\{0\} \times ]-h_1,0[$,
passes through $(w_1 - l/2 + \lambda_n , 0) \in ]w_1 -l ,
w_1[\times \{ 0 \}$, passes $p_n$ times through the line
$\{w_1-l\} \times ]0,h_2[$ and then passes through $A$.
So, $\gamma_n$ lies completely in $\C_1 \cup \C_2 \cup ]w_1-l,w_1[ \times \{0\}$.\\

\begin{figure}[!h]
\begin{center}
\psfrag{ln}{\small $\lambda_n$}
\includegraphics[width=.6 \linewidth]{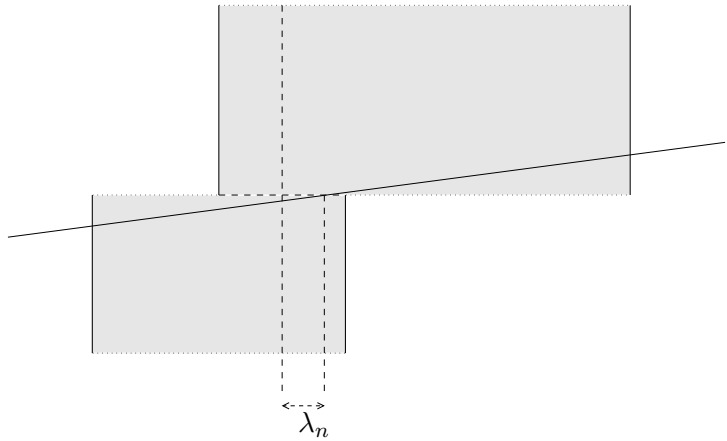}
\caption{\label{zoom-lemme.eps} Zoom on Figure
\ref{unfold-lemme.eps}.}
\end{center}
\end{figure}

Now, we assume by contradiction that there is a point $B(x,y)$ in
$\S$ distinct from $O$ and $A$ such that infinitely many
$\gamma_n$ pass through $B$. Hence, there is a subsequence such
that for all $n$ in $\N$, $\gamma_{i_n}$ passes through $B$. There
are two cases to consider:\\

First case: $y \in ]-1,0]$. By looking at the unfolded version of
the trajectory (Figure \ref{unfold-lemme.eps}), we see that, if
$k_{i_n}$ denotes the number of times that $\gamma_{i_n}$ pass
through the line $\{0\} \times ]-h_1,0[$ before hitting $B$, then
(by calculating the slope of $\gamma_{i_n}$ from $O$ to $B$)
$$\frac{y+1}{( x + k_{i_n} w_1 )-(w_1-l/2)} = \frac{2}{q_n w_1  +  p_n w_2 }.$$

So, $$x - w_1 + l/2 = ( q_{i_n} w_1 +  p_{i_n} w_2) (y+1) / 2 -
k_{i_n} w_1 .$$

In particular,
\begin{eqnarray*}
x - w_1 + l/2 & = & ( q_{i_0} w_1 +  p_{i_0} w_2) (y+1) / 2 - k_{i_0} w_1\\
              & = & ( q_{i_1} w_1 +  p_{i_1} w_2) (y+1) / 2 - k_{i_1} w_1.
\end{eqnarray*}

Hence, $$(p_{i_1}-p_{i_0}) w_2 + (q_{i_1}-q_{i_0}) w_1 = \frac{2
w_1 }{y+1} (k_{i_1}-k_{i_0}) \neq 0 .$$

So, $\frac{ 2 w_1}{y+1}$ can be written as $r w_2+s w_1$ where $r$
and $s$ are rational numbers.

Now, if $n\geq 1$, we still have $$(p_{i_n}-p_{i_0}) w_2 +
(q_{i_n}-q_{i_0}) w_1 = (r w_2 + s w_1) (k_{i_n}-k_{i_0}) .$$

Because $(w_2,w_1)$ is free over $\Q$, we have
\begin{itemize}
\item $(p_{i_n}-p_{i_0}) = r (k_{i_n}-k_{i_0})$ \item
$(q_{i_n}-q_{i_0}) = s (k_{i_n}-k_{i_0}) \neq 0$ (remember that
$q_n$ is strictly increasing)
\end{itemize}

Thus, by dividing, $$\frac{r}{s} =
\frac{p_{i_n}-p_{i_0}}{q_{i_n}-q_{i_0}} = \frac{p_{i_n}}{q_{i_n}}
(1-\frac{p_{i_0}}{p_{i_n}}) (\frac{1}{1-\frac{q_{i_0}}{q_{i_n}}})
\cge \frac{w_1}{w_2} \in \R \setminus \Q $$

leading to a contradiction.

For the second case, if $y \in [0,1[$, it is exactly the same
(just reverse Figure \ref{unfold-lemme.eps}).

Thus, $(O,A)$ is not finitely blockable, and $\S$ fails the finite
blocking property.

\hfill $\square$ \\

\clearpage

In the proof of Lemma \ref{lem2cyln}, we have chosen the points
$O\in \C_1$ and $A\in \C_2$ in such a way that the calculus was
easy. In fact, {\em almost all} pair $(O,A)\in \C_1 \times \C_2$
is not finitely blockable. More precisely, we have the

\begin{prop} \label{prop-local-rafinement}

In the system of coordinates given in Figure \ref{2cylinder.eps},
if $\alpha \defi w_1 / w_2 \in \R \setminus \Q$, if $A_i(x_i,y_i)
\in \C_i$ ($i \in \{1,2\}$), if $\beta \defi - y_1 / y_2$ is
uncommensurable to $\alpha$, then $(A_1,A_2)$ is not finitely
blockable.


\end{prop}

{\bf Proof:} it is very similar to the proof of Lemma
\ref{lem2cyln}; the introduction of many parameters just makes it
heavier to read.

Let $R \defi -y_2 l + y_2 w_1 -y_2 x_1 + x_2 y_1 - w_1 y_1 +y_1 l
$ (it is just an offset).

Let $$I \defi  ] \frac{-R}{w_2 y_2} ,\frac{l(y_2 - y_1) - R}{w_2
y_2} [.$$

$I$ is a non empty open set since $l(y_2 - y_1)>0$.

Since $\frac{\alpha}{\beta}$ is a positive irrational number, $\N
- \frac{\alpha}{\beta} \N$ is dense in $\R$ so there exists two
positive integer sequences $(p_n)_{n \in \N}$ and $(q_n)_{n \in
\N}$ such that:

\begin{itemize}
  \item $q_n$ is strictly increasing
  \item $\alpha q_n - \beta p_n \in I$
\end{itemize}


Let $$\lambda_n \defi \frac{y_2 w_1 q_n + y_1 w_2 p_n + R}{y_2 -
y_1} \in ]0,l[ \ \ \ (n \in \N).$$

For $n \in \N$, let $\gamma_n$ be the geodesic starting form $A_1$
with slope
$$\frac{y_2-y_1}{ w_1 q_n +  w_2 p_n + x_2 -x_1}
=\frac{-y_1}{ w_1 q_n + w_1 -x_1 -l +\lambda_n}=\frac{y_2}{ w_2
p_n - w_1 + x_2 +l  -\lambda_n}.$$

\begin{figure}[!h]
\begin{center}
\psfrag{totalslope}{\small $ w_1 q_n +  w_2 p_n + x_2 -x_1$}
\psfrag{slope1}{\small $ w_1 q_n + w_1 -x_1 -l +\lambda_n$}
\psfrag{slope2}{\small $w_2 p_n - w_1 + x_2 +l  -\lambda_n$}
\psfrag{ln}{\small $\lambda_n$} \psfrag{x1}{\small $x_1$}
\psfrag{y1}{\small $-y_1$} \psfrag{X2}{\small $w_1 + w_2 - l -
x_2$} \psfrag{Y2}{\small $y_2$}
\includegraphics[width=.9 \linewidth]{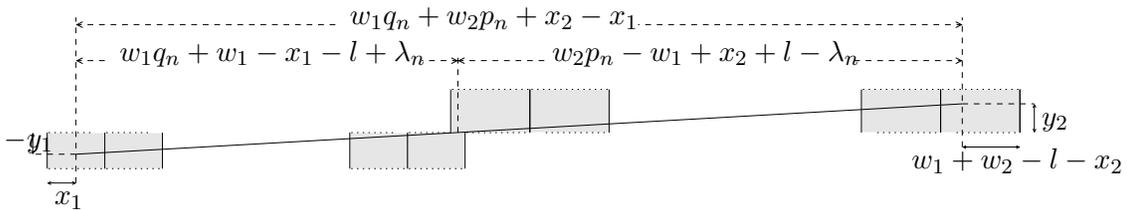}
\caption{\label{cunfold} The unfolding procedure.}
\end{center}
\end{figure}

So, we can check that $\gamma_n$ passes $q_n$ times through the
line $\{0\} \times ]y_1,0[$, passes through $(x_1 - l + \lambda_n
,0) \in ]x_1 -l , x_1[\times \{ 0 \}$, passes $p_n$ times through
the line $\{x_1-l\} \times ]0,y_1[$ and passes through $A_2$. Stop
$\gamma_n$ here (see Figure \ref{cunfold}).
So, $\gamma_n$ lies completely in $\C_1 \cup \C_2 \cup ]x_1-l,x_1[ \times \{0\}$.\\

Now, assume by contradiction that there is a point $B(x,y)$ in
$\S$ distinct from $A_1$ and $A_2$ such that infinitely many
$\gamma_n$ pass through $B$. Let $i \in \extr$ be an extraction
such that, for all $n$ in $\N$, $\gamma_{i_n}$ pass through $B$.
There are two cases to consider:

First case: $y \in ]y_1,0]$. By looking at the unfolded version of
the trajectory (Figure \ref{cunfold}), we see that
if $k_{i_n}$ denotes the number of times that $\gamma_{i_n}$ pass
through the line $\{0\} \times ]y_1,0[$ before hitting $B$, then
$$\frac{y_1-y}{x_1-(x+k_{i_n} w_1)} = \frac{y_2-y_1}{ w_1 q_n +  w_2 p_n + x_2 -x_1}.$$

So, $x-x_1 = ( w_1 q_{i_n} +  w_2 p_{i_n} + x_2 -x_1) (y-y_1) /
(y_2-y_1 ) - k_{i_n} w_1 $.



In particular,
\begin{eqnarray*}
x -  x_1 & = & ( w_1 q_{i_0} +  w_2 p_{i_0} + x_2 -x_1) (y-y_1) / (y_2-y_1 ) - k_{i_0} w_1 \\
         & = & ( w_1 q_{i_1} +  w_2 p_{i_1} + x_2 -x_1) (y-y_1) / (y_2-y_1 ) - k_{i_1} w_1.
\end{eqnarray*}

Hence, $(p_{i_1}-p_{i_0}) + (q_{i_1}-q_{i_0}) \alpha = \alpha
\frac{ y_2-y_1 }{y-y_1} (k_{i_1}-k_{i_0}) \neq 0$.

So, $\alpha \frac{ y_2-y_1 }{y-y_1}$ can be written as $r+s\alpha$
where $r$ and $s$ are rational numbers.

Now, if $n\geq 1$, we still have $(p_{i_n}-p_{i_0}) +
(q_{i_n}-q_{i_0}) \alpha = (r+s\alpha) (k_{i_n}-k_{i_0})$.

Because $(1,\alpha)$ is free over $\Q$, we have
\begin{itemize}
\item $(p_{i_n}-p_{i_0}) = r (k_{i_n}-k_{i_0})$ \item
$(q_{i_n}-q_{i_0}) = s (k_{i_n}-k_{i_0}) \neq 0$ (remember that
$q_n$ is strictly increasing)
\end{itemize}

Thus, by dividing, $$\frac{r}{s} =
\frac{p_{i_n}-p_{i_0}}{q_{i_n}-q_{i_0}} = \frac{p_{i_n}}{q_{i_n}}
(1-\frac{p_{i_0}}{p_{i_n}}) (\frac{1}{1-\frac{q_{i_0}}{q_{i_n}}})
\cge \frac{\alpha}{\beta} \in \R \setminus \Q $$

leading to a contradiction.

For the second case, if $y \in [0,y_2[$, it is exactly the same
(just reverse Figure \ref{cunfold}).

Thus, $(A_1,A_2)$ is not finitely blockable.

\hfill $\square$ \\

\clearpage

\mathversion{bold}
\section{Finite blocking property in the regular polygons}
\mathversion{normal}

\begin{thm}\label{th-regular}
Let $n\geq 3$ be an integer. The following assertions are
equivalent:
\begin{itemize}
\item the regular $n$-gon has the finite blocking property. \item
the right-angled triangle $\T_{n}$ with an angle equal to  $\pi/n$
has the finite blocking property. \item $n\in\{3,4,6\}$.
\end{itemize}
\end{thm}

{\bf Proof:} according to propositions \ref{prop-depliage} and
\ref{bound-almost}, it suffices to prove that the second assertion
implies the third one.

\begin{figure}[!h]
\begin{center}
\psfrag{s0}{$s_0$} \psfrag{s1}{$s_1$} \psfrag{s2}{$s_2$}
\psfrag{s-n/2}{$s_{\lfloor n/2 \rfloor}$}
\psfrag{s+n/2}{$s_{\lceil n/2 \rceil}$} \psfrag{sn-1}{$s_{n-1}$}
\psfrag{sn-2}{$s_{n-2}$} \psfrag{s'0}{$s'_0$} \psfrag{s'1}{$s'_1$}
\psfrag{s'2}{$s'_2$} \psfrag{s'n-1}{$s'_{n-1}$}
\psfrag{s'n-2}{$s'_{n-2}$}
\includegraphics[scale=0.5]{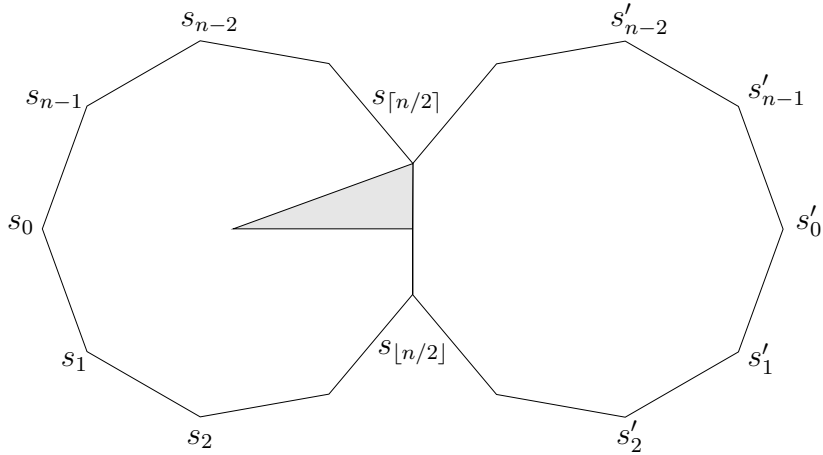}
\caption{\label{nona.eps} The translation surface associated with
the right-angled triangle with an angle equal to $\pi/9$ (identify
the opposite sides).}
\end{center}
\end{figure}

Suppose first that $n$ is odd. $ZK(\T_{n})$ can be described as
two regular $n$-gons $P$ and $P'$ symmetric one to each other,
with identifications along the sides (see Figure \ref{nona.eps}).

We number the vertices of $P$ (resp. $P'$) counterclockwise (resp.
clockwise) from $s_0$ to $s_{n-1}$ (resp. from $s'_0$ to
$s'_{n-1}$) (see Figure \ref{nona.eps}).


\begin{figure}[!h]
\begin{center}
\psfrag{cut}{cut} \psfrag{d2}{$d_2$} \psfrag{d4}{$d_4$}
\psfrag{paste}{paste} \psfrag{rotate}{rotate}
\psfrag{stretch}{stretch}
\includegraphics[width=.9 \linewidth]{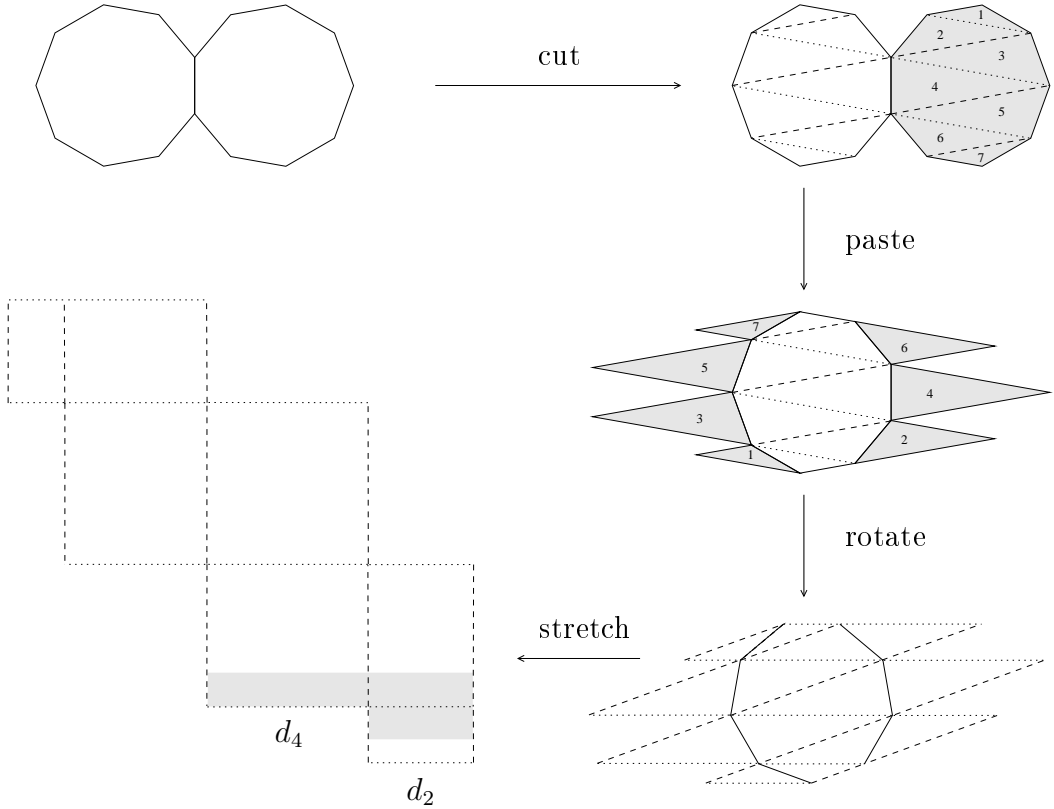}
\caption{\label{nona-fusion.eps} From the surface associated to
$\T_{n}$ to a more exploitable one in the same orbit under
$GL(2,\R)$ ($n$ odd).}
\end{center}
\end{figure}

Now, look at Figure \ref{nona-fusion.eps}.

In $P\cup P'$, $(s_0,s_{\lceil n/2 \rceil})$ is parallel to
$(s_1,s_{\lceil n/2 \rceil}-1)$, to $(s_2,s_{\lceil n/2
\rceil}-2)$, \dots, to $(s_{n-1},s_{\lceil n/2 \rceil}+1)$, to
$(s_{n-2},s_{\lceil n/2 \rceil}+2)$, \dots, and by axial symmetry
with respect to $(s_{\lfloor n/2 \rfloor},s_{\lceil n/2 \rceil})$,
it is parallel to $(s'_0,s'_{\lceil n/2 \rceil})$, to
$(s'_1,s'_{\lceil n/2 \rceil}-1)$, to $(s'_2,s'_{\lceil n/2
\rceil}-2)$, \dots, to $(s'_{n-1},s'_{\lceil n/2 \rceil}+1)$, to
$(s'_{n-2},s'_{\lceil n/2 \rceil}+2)$, \dots.

Let us call this common direction the dashed one.

Do the same for the direction $(s_0,s_{\lfloor n/2 \rfloor})$ and
call it the dotted one.

This leads to a triangulation of the surface, each triangle having
an edge dashed, an edge dotted and an edge that is an edge of $P$
or $P'$.

Now, take $P$ as the base and glue all the triangles that lie in
$P'$ to the ones that lie in $P$, thanks to the identification
between the edges of $P$ and the edges of $P'$ (this can be done
in an unique way) (this is the cut and paste operation).

We have now a new representation of the surface associated to
$\T_{n}$, by a planar polygon with identifications along its sides
that are only dashed or dotted sides.

Since those two directions are not parallel, it is easy to find an
element $A$ of $GL(2,\R)$ that put the dotted direction to the
horizontal, the dashed one to the vertical and that preserves the
lengths in those two directions (this is the rotate and stretch
operation).

The surface $\S$ obtained by the action of $A$ to the surface
associated to $\T_{n}$ is more exploitable for our purpose.

For $1\leq i \leq n-1$, we denote $d_i \defi || s_0 - s_i ||_2$.

We can remark that exactly one edge of $P$ is dashed. Starting
from this edge in $\S$, it is easy to recognize the shape of the
Figure \ref{2cylinder.eps} with $w_1=l=d_2$, $w_2=d_2+d_4$,
$h_1=d_1$ and $h_2= d_3$.

We have $w_2 / w_1 = 1+d_4 / d_2 = 1+ 2\sin(4\pi / n ) /
2\sin(2\pi / n ) = 1+ 2\cos(2\pi / n )$ which is irrational if
$n\geq 5$. Hence, $\S$ and therefore $\T_{n}$ lacks the finite
blocking property if $n\neq 3$.\\ \\



\begin{figure}[!h]
\begin{center}
\psfrag{cut}{cut} \psfrag{d2}{$d_2$} \psfrag{d4}{$d_4$}
\psfrag{paste}{paste} \psfrag{rotate}{rotate}
\psfrag{stretch}{stretch}
\includegraphics[width=.9 \linewidth]{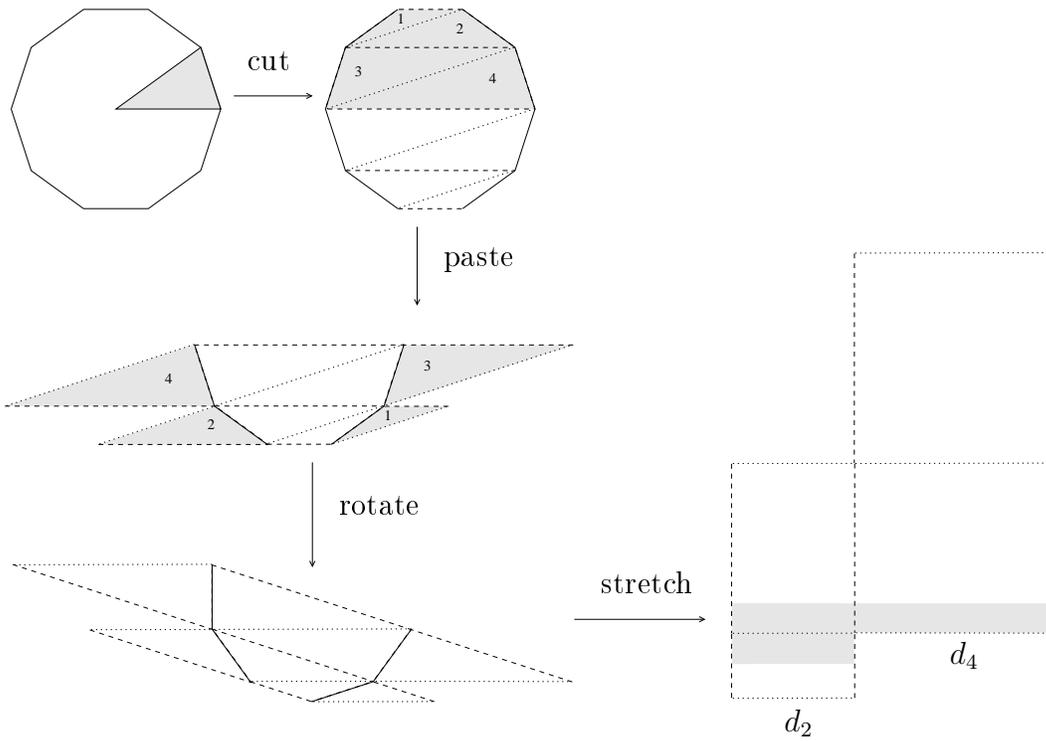}
\caption{\label{deca-fusion.eps} The even case ($n=10$).}
\end{center}
\end{figure}

For the even case, the translation surface $ZK(\T_{n})$ is
constituted by only one regular $n$-gon with opposite sides
identified. The construction is very similar, but the role of $P$
is played by the lower-half of the polygon, and the role of $P'$
is played by its upper-half (see Figure \ref{deca-fusion.eps}). In
this case, $w_2 / w_1 = 1+ 2\cos(2\pi / n )$ is irrational if
$n\notin \{4,6\}$.

\hfill $\square$ \\


Nevertheless, we can notice that the situation is {\em not
homogenous} among different pairs of points:

\begin{prop}
For every regular $n$-gon ($n$ even), there exists a finite set of
points that blocks every billiard path from the center $O$ to
itself.
\end{prop}


{\bf Proof:} the set $\B$ consisting in the centers of the edges
and the vertices is a (finite) blocking configuration. Indeed,
suppose by contradiction that $\gamma$ is a billiard path from $O$
to $O$ that does not meet $\B$: it can be folded into a billiard path
in the triangle $\T_n$ from the vertex with angle $\pi / n$ to
himself. This contradicts  \cite{To}, Lemma 4.1. page 871 (a nice
argument on angles, measured modulo $2\pi / n$ makes this fact
impossible).

\hfill $\square$ \\

\clearpage



\begin{figure}[!h]
\begin{center}
\includegraphics[width=.9 \linewidth]{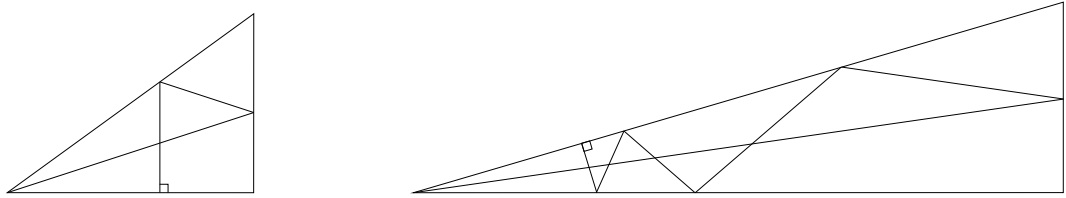}
\caption{\label{noblock.eps} By starting with an angle $\pi/2n$,
we are coming back in $O$ ($n\in \{5,11\} \subset 2\N + 3$).}
\end{center}
\end{figure}

This proposition is not so clear if $n$ is odd, since a billiard
path $\gamma$ in $\T_n$ starting from $O$ with angle $\pi/2n$ is
coming back to $O$ after $n$ bounces (see Figure
\ref{noblock.eps}).

It would be interesting to find a rational polygonal billiard
table $\P$ with the property that no pair of points in $\P$ can be
finitely blocked.

\clearpage

\section{Finite blocking property on Veech surfaces} \label{2generalisation}

\begin{prop}
Suppose that a translation surface $\S$ is decomposable into
commensurable parallel cylinders, in at least two different
directions.

Then $\S$ has the finite blocking property if and only if $\S$ is
a torus covering, branched over only one point.
\end{prop}

{\bf Proof:} Let $(\C_i)_{i=0}^n$ be a decomposition of $\S$ into
parallel commensurable cylinders, with heights $h_i$ and weights
$w_i$ ($i\leq n$).

Starting from $\C_0$, by applying Lemma \ref{lem2cyln} step by
step, we can see that all the $w_i$'s are commensurable (recall
that $\S$ is assumed to be connected).

Since, for all $i\leq n$, $w_i / h_i$ is a rational number, we can
deduce that all the $h_i$'s are commensurable. So, there exists
$h>0$ and $(k_i)_{i=0}^n \in \N^{n+1}$ such that for $i\leq n$,
$h_i=k_i h$. Then, each $\C_i$ is decomposed into $k_i$
subcylinders $(\C_{i,j})_{j=1}^{k_i}$ of height h ($\C_{i,j}$ is
the image of $\R / w_i \Z \times ](j-1) h , j h[$). Note that the
singularities of $\S$ lie in the sides of the $\C_{i,j}$'s.

By hypothesis, we have the same kind of decomposition of $\S$ into
parallel subcylinders $(\C'_{i,j})_{i\leq n',j\leq k'_i}$ of
height $h'$ in another direction, with the property that each
singularity of $\S$ lies in a side of $\C'_{i,j}$.

$$(\P_i)_{i\leq l} \defi (\C_{i,j})_{i\leq n,j\leq k_i}  \bigvee  (\C'_{i,j})_{i\leq n',j\leq k'_i}$$
is a decomposition of $\S$ into parallel isometric parallelograms
glued edge to edge.

This leads to a covering from $\S$ to $\P_0$ whose opposite edges
are identified, i.e. from $\S$ to a torus. Note that all the
singularities of $\S$ lie in a vertex of some $\P_i$, so they are
sent to a common point in the torus (the image of a vertex of
$\P_0$).

\hfill $\square$ \\

\begin{thm}\label{th-Veech}
A Veech surface has the finite blocking property if and only if it
is a torus covering, branched over only one point.
\end{thm}

{\bf Proof:} If the surface is a torus, there is nothing to prove
(both statements are true). Otherwise, the genus of the surface is
greater than two, and then the surface has at least one
singularity and therefore many saddle connection directions. In
each saddle connection direction, a Veech surface admits a
decomposition into commensurable cylinders.

\hfill $\square$ \\


\clearpage

\section{Further results}

The aim of this section is to see how Lemma \ref{lem2cyln} can be
used in different contexts.

\begin{figure}[!h]
\begin{center}
\psfrag{w1}{$w_1$} \psfrag{a}{$\alpha$} \psfrag{w2}{$w_2$}
\psfrag{1}{$1$}
\includegraphics[scale = 0.9]{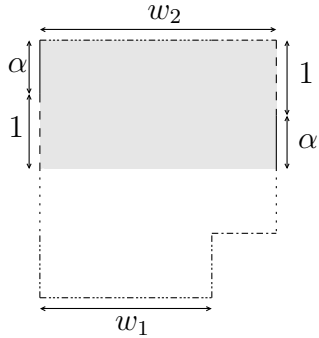}
\caption{\label{localcylinder.eps} The surface cannot be fully
decomposed into cylinders in the horizontal direction (the grey
zone is a minimal component for the horizontal flow); identify the
lines with the same style, $\alpha \in \R \setminus \Q$, $w_1/ w_2
\in \R \setminus \Q$.}
\end{center}
\end{figure}

First, recall that Lemma \ref{lem2cyln} can be applied in the
non-Veech context since:
\begin{itemize}
\item $\C_1$ and $\C_2$ do not need to be commensurable (there is
no condition on the heights) \item the result is {\em local}: if
the configuration of Figure \ref{2cylinder.eps} appears somewhere
in a surface $\S$ with $w_1/ w_2 \in \R \setminus \Q$, then $\S$
lacks the finite blocking property (see example on Figure
\ref{localcylinder.eps}).
\end{itemize}

\subsection{L-shaped surfaces}

\begin{figure}[!h]
\begin{center}
\psfrag{1}{$1$} \psfrag{a}{$a$} \psfrag{b}{$b$}
\includegraphics[scale=0.5]{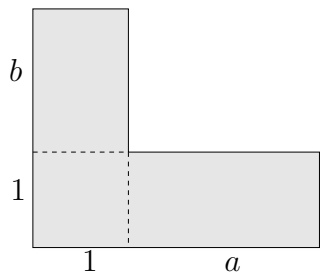}
\caption{\label{L-shaped.eps} The L-shaped translation surface
$L(a,b)$ (identify the opposite sides by translation).}
\end{center}
\end{figure}

Let $\S$ be a \mbox{L-shaped} translation surface; it is in the
same $GL(2,\R)$-orbit than some $L(a,b)$ (see Figure
\ref{L-shaped.eps}). Lemma \ref{lem2cyln} allows us to decide
wether $\S$ has the finite blocking property (we do not need $\S$
to be a Veech surface (such surfaces were characterized in
\cite{Mc} and \cite{Cal})).

\begin{prop} \label{L-prop}
Let $a$ and $b$ be two positive real numbers. Then $L(a,b)$ has
the finite blocking property if and only if $(a,b)\in \Q^2$.
\end{prop}

{\bf Proof:} If $L(a,b)$ has the finite blocking property,
applying Lemma \ref{lem2cyln} in both horizontal and vertical
direction leads to $(1+a)/1\in\Q$ and $(1+b)/1\in\Q$, so $(a,b)\in
\Q^2$.

\hfill $\square$ \\

\subsection{Irrational billiards}

The construction of Zemljakov and Katok is also possible when the
angles of a polygon $\P$ are not rational multiples of $\pi$; in
this case, the group $\Gamma$ is infinite and the surface is not
compact. However, periodic trajectories can appear; we can even
meet the situation of lemma \ref{lem2cyln}:

\begin{prop}
There exists a non rational polygonal billiard that fails the
finite blocking property.
\end{prop}

\begin{figure}[!h]
\begin{center}
\psfrag{1}{$1$} \psfrag{a}{$\alpha$} \psfrag{t}{$\theta$}
\includegraphics[scale=0.6]{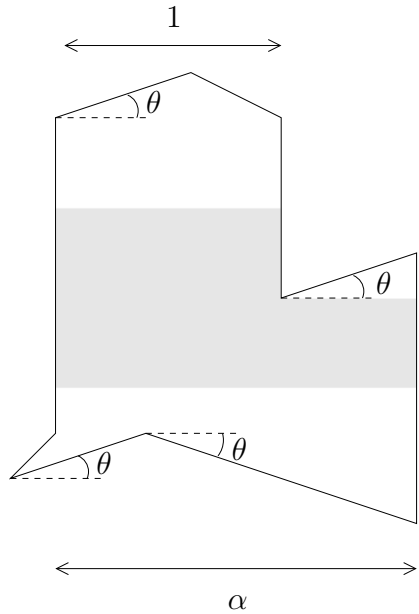}
\caption{\label{irrational} A non rational billiard whose
associated surface contains the configuration of Lemma
\ref{lem2cyln}; $\theta \in \R \setminus \Q \pi$, $\alpha \in \R
\setminus \Q$.}
\end{center}
\end{figure}

{\bf Proof:} consider the billiard drawn in Figure
\ref{irrational} and apply Lemma \ref{lem2cyln} on the pair of
cylinders defined by the grey zone.

\hfill $\square$ \\

\subsection{A density result}

A singularity $\sigma \in \Sigma$ has a conical angle of the form
$2k\pi$, with $k\geq 1$; we say that $\sigma$ is of {\em
multiplicity} $k-1$. If $1\leq k_1\leq k_2 \leq \dots \leq k_n$ is
a sequence of integers whose sum is even, we denote $\H(k_1, k_2,
\dots, k_n)$ the {\em stratum} of translation surfaces with
exactly $n$ singularities whose multiplicities are
$k_1, k_2, \dots, k_n$. A translation surface in $\H(k_1, k_2,
\dots, k_n)$ has genus $g=1+(k_1+k_2+ \dots +k_n)/2$.

Each stratum carries a natural topology that is for example
defined in \cite{Ko}.

\begin{thm}\label{th-density}
In genus $g\geq 2$, the set of translation surfaces that fail the
finite blocking property is dense in every stratum.
\end{thm}

{\bf Sketch of the proof:} The proof of this requires some
material that is too long to describe here (like the precise
definition of the topology on such a stratum (for this we have to
see each translation surface as a Riemann surface with an abelian
differential)). It suffices to prove that the translation surfaces
that satisfy the hypothesis of Lemma \ref{lem2cyln} are dense in
each stratum. For this, we begin to prove that translation
surfaces that admit a cylinder decomposition in the horizontal
direction, with at least two non homologous horizontal cylinders,
are dense (see \cite{EO}, \cite{KZ}, \cite{Zo}). Using the local
coordinates given by the period map, we have the possibility to
perturb the perimeters of two such ``consecutive'' non homologous
cylinders in order to let them non commensurable. We postpone the
precise proof for a further paper (see \cite{Mo-future}).

\hfill $\square$ \\

In fact, we will prove in \cite{Mo-future} that the translation
surfaces that fail the finite blocking property is of full measure
in each stratum. Moreover, we will prove that finite blocking
property implies complete periodicity; we will also give the
classification of the surfaces that have the finite blocking
property in genus $2$.\\ \\


Of course,
the notion of finite blocking property and the results presented in this paper can be translated in the vocabulary of quadratic differentials.

\clearpage

\section*{Conclusion}

One can define a stronger property: a planar polygonal billiard or
a translation surface $\P$ is said to have the {\em bounded
blocking property} if the number of blocking points can be chosen
independently of the pair $(O,A)$. Does it exist polygonal
billiard
tables with the finite blocking property but without the bounded blocking property?\\

Is it true that for general translation surfaces, the fact of
being a torus covering is a necessary and sufficient condition to
have the finite blocking property (resp. bounded blocking
property)? Is it true that for rational billiards, the fact of
being almost integrable is a necessary and sufficient condition to have the finite blocking property
(resp. bounded blocking property)?\\


Note that for piecewise smooth billiard tables, the study of the
finite blocking property seems very difficult,
since for an ellipse, even a countable number of points do not suffice to block every path from a focus to the other one.\\

\end{document}